\documentclass[11pt]{amsart}
\usepackage{amssymb,amsmath,latexsym,verbatim,epsf,enumerate,graphicx}

\def\Z{\mathbb{Z}}
\def\N{\mathbb{N}}
\def\R{\mathbb{R}}
\def\C{\mathcal{C}}
\def\D{\mathcal{D}}
\def\f{f}
\def\b{{\bf b}}
\def\m{{\bf m}}
\def\v{{\bf v}}
\def\0{{\bf 0}}
\def\1{{\bf 1}}
\def\s{s}
\def\th{^{\text{th}}}
\newcommand\an[1]{\angle \left( {#1} \right)}
\renewcommand\sp[2]{\left( {#1} , {#2} \right)}
\newcommand\fl[1]{\left\lfloor {#1} \right\rfloor}
\newcommand\cl[1]{\left\lceil {#1} \right\rceil}

\newtheorem{theorem}{Theorem}
\newtheorem{definition}{Definition}

\newtheorem{conjecture}{Conjecture}
\newtheorem{lemma}[theorem]{Lemma}

\numberwithin{equation}{section}

\title[Higher-dimensional Dedekind sums and their bounds]{Higher-dimensional Dedekind sums and their bounds arising from the discrete diagonal of the $n$-cube}

\author{Matthias Beck}
\address{Department of Mathematical Sciences\\
        Binghamton University (SUNY)\\
        Binghamton, NY 13902-6000\\
        USA}
\email{matthias@math.binghamton.edu}

\author{Sinai Robins}
\address{Department of Mathematics\\
        Temple University\\
        Philadelphia, PA 19122\\
        USA}
\email{srobins@math.temple.edu}

\author{Shelemyahu Zacks}
\address{Department of Mathematical Sciences\\
        Binghamton University (SUNY)\\
        Binghamton, NY 13902-6000\\
        USA}
\email{shelly@math.binghamton.edu}

\thanks{The second author would like to thank the generous support of NSA Young Investigator grant MSPR-OOY-196}

\begin{document}
\setlength{\parindent}{0pt}
\maketitle

\small
\begin{quote}
{\it Counting pairs is the oldest trick in combinatorics... Every time we count pairs, we learn something from it.} \\
Gil Kalai
\end{quote}
\normalsize

\abstract Higher-dimensional Dedekind sums are defined as a
generalization of a recent 1-dimensional probability model of
Dilcher and Girstmair to a $d$-dimensional cube. The analysis of
the frequency distribution of marked lattice points leads to new
formulae in certain special cases, and to new bounds for the
classical Dedekind sums. Upper bounds for the generalized Dedekind
sums are defined in terms of 1-dimensional moments. In the
classical two-dimensional case, the ratio of these sums to their
upper bounds are cosines of angles between certain vectors of
n-dimensional cones, conjectured to have a largest spacial angle
of $\pi/6$.
\endabstract

\setlength{\parskip}{0.4cm}
\bibliographystyle{amsplain}


\section{Introduction}

Historically, Dedekind sums first appeared in Dedekind's transformation law of his $\eta$-function \cite{dedekind}.
Dedekind sums have since become an integral part of combinatorial geometry (lattice point enumeration \cite{mordell}),
algebraic number theory (class number formulae \cite{meyerdedekind}), topology (signature defects of manifolds \cite{hirzebruchzagier}),
and algorithmic complexity (pseudo random number generators \cite{knuth}).
We begin by defining the classical Dedekind sum, whose basic ingredient is the sawtooth function
  \[ ((x)) = \left\{ \begin{array}{cl} \{ x \} - \frac{ 1 }{ 2 } & \mbox{ if } x \not\in \Z \\
                                       0                         & \mbox{ if } x \in \Z \ . \end{array} \right. \]
Here $ \{ x \} = x - \fl x $ denotes the fractional part of a real
number $x$.

For any two positive integers $a$ and $b$, we define the \emph{classical Dedekind sum} as
  \begin{equation} \s (a,b) = \sum_{ k \text{ {\rm mod} } b } \left( \left( \frac{ ka }{ b } \right) \right) \left( \left( \frac{ k }{ b } \right) \right) \ . \end{equation}
Here the sum is over a complete residue system modulo $b$.

The classic introduction to the arithmetic properties of the Dedekind sum is \cite{grosswald}.
The Dedekind sums have recently been cast in a new light as essentially the second moments of
an appealing probability model introduced by Dilcher and Girstmair \cite{dilchergirstmair}.
They divide an interval of length $a$ into $b$ equal subintervals (``boxes") and
count the number of integers in each subinterval.

We generalize their approach by considering a $d$-dimensional cube ($d \geq 2$) of side
length $a \in \N$. Along the $d$-dimensional main diagonal we mark the points with
integer coordinates. The cube is now partitioned to $b_1 b_2 \cdots b_d$ rectangular
parallelopipeds by dividing the $j\th$ side of the cube to $b_j$ equal length intervals.
Each parallelopiped is given the coordinates $(j_1, j_2, \dots, j_d)$ where $j_k = 1, \dots, b_k \
(k=1, \dots, d)$. Let $f_{a; b_1, \dots, b_d} (j_1, \dots, j_d)$ denote the number (frequency) of marked points
along the main diagonal which belong to the $(j_1, j_2, \dots, j_d)$-parallelopiped. The
generalized Dedekind sums under consideration are
  \begin{equation} S_d (a; \b) = S_d (a; b_1, b_2, \dots, b_d) = \frac 1 a \sum_{k_1=0}^{b_1-1} \cdots \sum_{k_d=0}^{b_d-1} k_1 \cdots k_d \, \f_{a; b_1, \dots, b_d} \left( k_1, \dots, k_d \right)  , \end{equation}
a mixed moment for the $f_{a; b_1, \dots, b_d}$ distribution.
 In Section \ref{secondintro} we present basic definition and facts.
We show that
  \begin{equation} S_d (a; \b) = \frac 1 a \sum_{ m=0 }^{a-1} \fl{ \frac{ m b_1 }{ a } } \cdots \fl{ \frac{ m b_d }{ a } } \ , \end{equation}
For the case $d=1$ we define the $k\th$ moment as
  \begin{equation} M_k(a;b) = \frac 1 a \sum_{m=0}^{a-1} \left\lfloor \frac{ mb }{ a } \right\rfloor^k \ . \end{equation}
These moments are used in Sections \ref{d=2bounds} and \ref{boundshigherdim} to provide,
via the Cauchy-Schwartz inequality, upper bounds for $S_d (a; \b)$.

The second moment $M_2(a;b)$ is of special importance since it is connected with the classical
Dedekind sum $\s(a,b)$ according to the formula (see Section \ref{secondintro})
  \begin{equation}\label{M2} M_2(a;b) = \frac{ (b^2 + 1) (a-1) (2a-1) }{6a^2} - \frac{ (a-1)b }{ 2a } - \frac{ 2b }{a} \, \s(b,a) \ . \end{equation}
By analyzing the structure of the univariate and bivariate
frequency distributions (Section \ref{2dimsect}), we derive in
Section \ref{dilchersect} new formulae for $M_2(a;b)$, in some
special cases, and provide several types of lower and upper
bounds. In Section \ref{d=2bounds} we analyze the ratio of
$S_2(a;b_1,b_2)$ to its upper bound $\left( M_2(a;b_1) M_2(a;b_2)
\right)^{1/2}$. All these ratios $R_2(a;b_1,b_2)$ are empirically
found to be greater or equal to $R_2(5;2,3) = \sqrt 3 /2$, leading
to the following conjecture (see section $5.1$ below):

{\bf Conjecture 1}.  For all $a, b, c \geq 3$,
\begin{equation} R_2 (a;b,c) \geq R_2 (5;2,3) = \frac{ \sqrt 3 }{ 2 } \ . \end{equation}

Geometrically, $S_2(a;b_1,b_2)$ is an inner product of the vectors
$\v_{b_1} = \left( \fl{ \frac{b_1}{a} },  \fl{ \frac{ 2 b_1 }{a} }, \dots, \fl{\frac{(a-1)b_1}{a}} \right)$
and $\v_{b_2} = \left( \fl{ \frac{b_2}{a} }, \fl{ \frac{ 2 b_2 }{a} }, \dots, \fl{ \frac{(a-1)b_2}{a}} \right)$
in $\R^{a-1}$, and $R_2(a;b_1,b_2)$ is
the cosine of the angle between these two vectors. It appears from both empirical and theoretical evidence that
all these vectors, for $a \geq 3, \, b_1, b_2 \geq 2$, are within
a cone with largest possible angle of $\cos^{-1} ( \sqrt 3 / 2 ) =
\pi/6$. In Section \ref{cauchy-schwsubsec} we have some general
results and observations on the functions $R_2(a;b_1,b_2)$. In
Section \ref{geomsubsec} we analyze the geometry of the vectors
$\v_{b} = \left(\fl{ \frac{b}{a} },  \fl{ \frac{ 2 b }{a} }, \dots, \fl{\frac{(a-1)b}{a}}
\right)$ and prove several lemmas, which lend further credence to
the validity of Conjecture \ref{firstconjecture}. Finally, in
Section \ref{boundshigherdim} we present higher-dimensional upper
bounds for $S_d (a; \b)$ in terms of the $r$'th moments $M_r(a;b)$, and prove
that $M_r (a; b)$ is log-convex in $r$.

\section{Generalizing the Dilcher-Girstmair model}\label{secondintro}

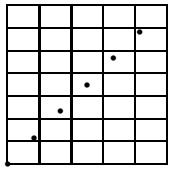
\begin{figure}[htb]
\begin{center}
\begin{picture}(60,60)
\put(0,0){\framebox(60,60)}
\put(12,0){\line(0,1){60}}
\put(24,0){\line(0,1){60}}
\put(36,0){\line(0,1){60}}
\put(48,0){\line(0,1){60}}
\put(0,8.75){\line(1,0){60}}
\put(0,17.14){\line(1,0){60}}
\put(0,25.71){\line(1,0){60}}
\put(0,34.29){\line(1,0){60}}
\put(0,42.86){\line(1,0){60}}
\put(0,51.43){\line(1,0){60}}
\put(0,0){\circle*{2}}
\put(10,10){\circle*{2}}
\put(20,20){\circle*{2}}
\put(30,30){\circle*{2}}
\put(40,40){\circle*{2}}
\put(50,50){\circle*{2}}
\end{picture}
\end{center}
\vspace{-15pt}
\caption{$a=6, \ b=5, \ c=7$} \label{example}
\end{figure}

We introduce the $d$-dimensional analog of the Dilcher-Girstmair model.
We begin gently with the two-dimensional
extension: given three positive integers $a$, $b$, and $c$,
divide one of the sides of the square $[0,a) \times [0,a)$ into $b$ parts of length $a/b$,
and the other side into $c$ parts of length $a/c$.
This division induces a grid (see Figure \ref{example} for an example).
We thus have $bc$ boxes of equal size.
We think of each box as half open: we count the left (excluding the highest point) and
bottom side (excluding the right-most point) as belonging to the box.
Let's mark each box by a pair of integers $(j,k)$ where $0 \leq j \leq b-1$
and $0 \leq k \leq c-1$. We will study the integer lattice points in the square; note that the
box $(j,k)$ contains the point $(m,n) \in \Z^2$ if and only if
  \begin{equation}\label{one} \frac{ ja }{ b } \leq m < \frac{ (j+1)a }{ b } \qquad \text{ and } \qquad \frac{ ka }{ c } \leq n < \frac{ (k+1)a }{ c } \ . \end{equation}
Equivalent to this condition is the following:
  \[ j \leq \frac{ mb }{ a } < j+1 \qquad \text{ and } \qquad k \leq \frac{ nc }{ a } < k+1 \ , \]
which can be rewritten in compact form using the greatest integer function $\lfloor x \rfloor$ (the greatest integer not exceeding $x$):
  \begin{equation}\label{two} j = \left\lfloor \frac{ mb }{ a } \right\rfloor \qquad \text{ and } \qquad k = \left\lfloor \frac{ nc }{ a } \right\rfloor \ . \end{equation}
We formalize the distribution of integer points within each of the $bc$ boxes as follows.

\begin{definition}
Let $\f_{a;b,c} (j,k)$ denote the number of marked lattice points on the diagonal in the box $(j,k)$.
\end{definition}
Notice that most of these frequencies are zero.
We can evaluate the following sum in two ways according to the equivalence of (\ref{one}) and (\ref{two}):
  \begin{equation}\label{twodimequ} \sum_{j=0}^{b-1} \sum_{k=0}^{c-1} j \, k \, \f_{a;b,c} (j,k) = \sum_{m=0}^{a-1} \left\lfloor \frac{ mb }{ a } \right\rfloor \left\lfloor \frac{ mc }{ a } \right\rfloor \ . \end{equation}
A special case of this is $b=c$, for which we make the following definition.

\begin{definition}
Let $f_{a;b} (j)$ denote the number of marked integers in the $j\th$ subinterval when we divide
$[0,a)$ into $b$ equal parts.
\end{definition}
We remark that $f_{a;b} (j) = f_{a;b,b} (j) $. When $b=c$ the double sums
reduce to the one-dimensional sums studied by Dilcher and Girstmair, i.e.
  \begin{equation}\label{m2equiv} \sum_{j=0}^{b-1} \sum_{k=0}^{b-1} j \, k \, \f_{a;b,b} (j,k) = \sum_{j=0}^{b-1} j^2 f_{a;b} (j) = \sum_{m=0}^{a-1} \left\lfloor \frac{ mb }{ a } \right\rfloor^2 \ . \end{equation}
The sum on the right hand is essentially a classical Dedekind sum:
If $a$ and $b$ are relatively prime,
  \begin{align*} \sum_{m=0}^{a-1} \left\lfloor \frac{ mb }{ a } \right\rfloor^2 &= \sum_{m=1}^{a-1} \left( \frac{ mb }{ a } - \left\{ \frac{ mb }{ a } \right\} \right)^2 \\
     &= \sum_{m=1}^{a-1} \left( \frac{ mb }{ a } \right)^2 - 2 \sum_{m=1}^{a-1} \frac{ mb }{ a } \left\{ \frac{ mb }{ a } \right\} + \sum_{m=1}^{a-1} \left\{ \frac{ mb }{ a } \right\}^2 \\
     &= \frac{ b^2 }{a^2} \sum_{m=1}^{a-1} m^2 - 2b \sum_{m=1}^{a-1} \left( \frac{ m }{ a } - \frac 1 2 \right) \left( \left\{ \frac{ mb }{ a } \right\} - \frac 1 2 \right) - 2b \sum_{m=1}^{a-1} \frac{ m }{ a } + \frac{2b(a-1)}{4} + \sum_{m=1}^{a-1} \left( \frac{ m }{ a } \right)^2 \\
     &= \frac{ (b^2 + 1) (a-1) (2a-1) }{6a}  - 2b \, \s(b,a) - \frac 1 2 b (a-1)
  \end{align*}
This and similar sums coming from the one-dimensional case will appear repeatedly in the
exposition that follows.

\begin{definition}
For any two positive integers $a$ and $b$, let
  \[ M_k(a;b) = \frac 1 a \sum_{m=0}^{a-1} \left\lfloor \frac{ mb }{ a } \right\rfloor^k \ . \]
\end{definition}
$M_k$ is the $k\th$ moment of the Dilcher-Girstmair probability distribution (see Section \ref{boundshigherdim}).
By (\ref{m2equiv}), the definition of $M_k$ is equivalent to
  \[ M_k(a;b) = \frac 1 a \sum_{j=0}^{b-1} j^k f_{a;b} (j) \ . \]
We just showed above that $M_2$ corresponds to the classical Dedekind sum $\s(a,b)$ as in (\ref{M2}).

The model that we described above extends naturally to higher dimensions.
Instead of considering a square, let's divide the $d$-dimensional cube $[0,a) \times \cdots \times [0,a)$
into $b_1 \cdots b_d$ equal boxes by a similar construction as above: now we divide
the first side into $b_1$ equal intervals, the next one into $b_2$ equal intervals, and
so on. Again we will count the number of marked integer lattice point on the main diagonal of this
cube, according to the box they are in. As above we will label each box, say by
$ \left( k_1, \dots, k_d \right) , \ 0 \leq k_j \leq b_j - 1$,
and we will denote the function counting the lattice
points in box $ \left( k_1, \dots, k_d \right) $ by
  \[ \f_{ a; b_1, \dots, b_d } \left( k_1, \dots, k_d \right) \ . \]
As before, an elementary counting-two-ways argument yields
  \[ \sum_{k_1=0}^{b_1-1} \cdots \sum_{k_d=0}^{b_d-1} k_1 \cdots k_d \, \f_{a; b_1, \dots, b_d} \left( k_1, \dots, k_d \right) = \sum_{m=0}^{a-1} \left\lfloor \frac{ mb_1 }{ a } \right\rfloor \cdots \left\lfloor \frac{ mb_d }{ a } \right\rfloor \ . \]
This naturally leads to the following definition.
\begin{definition} For positive integers $a, b_1, b_2, \dots, b_d$, we define
  \[ S_d (a; \b) = S_d (a; b_1, b_2, \dots, b_d) = \frac 1 a \sum_{ m=0 }^{a-1} \fl{ \frac{ m b_1 }{ a } } \cdots \fl{ \frac{ m b_d }{ a } } \ , \]
\end{definition}
This is a generalized Dedekind sum.
Our goal is to find relations for the sums $S_d (a; \b)$.


\section{The two-dimensional frequency distribution $\left\{ \f_{a;b,c} (j,k) \right\}$ and its marginal distributions}\label{2dimsect}

In this section we focus on the study of the distribution frequencies $\f_{a;b,c} (j,k)$
using the duality interpretation given by (\ref{twodimequ}).
It appears impossible to derive a closed formula for the number of diagonal lattice points that
belong to the $(j,k)$th rectangle, that is, $\f_{a;b,c} (j,k), \ j=0, \dots, b-1, \ k=0, \dots, c-1$.
We developed an algorithm, given in the appendix, for computing the values of $\f_{a;b,c} (j,k)$
and of the marginal frequencies
  \[ \f_{a;b} (j) = \sum_{k=0}^{c-1} \f_{a;b,c} (j,k) \qquad \text{ and } \qquad \f_{a;c} (k) = \sum_{j=0}^{b-1} \f_{a;b,c} (j,k) \ . \]

\emph{Example.}
In Table \ref{distribution50137} we present these
distributions for the case $a=50, b=13, c=7$. From this table we can immediately verify that
  \[ \sum_{j=0}^{12} \sum_{k=0}^6 k \, j \, f_{50;13,7} (j,k) = \sum_{m=0}^{49} \fl{ \frac{ 13 m }{ 50 } } \fl{ \frac{ 7 m }{ 50 } } = 1236 \ . \]
\begin{figure}[htb]
\begin{center}
\begin{tabular}{c|ccccccc|c}
$j \backslash k$ & 0 & 1 & 2 & 3 & 4 & 5 & 6 & $\f_b^a (j)$ \\ \hline
0                & 4 & 0 & 0 & 0 & 0 & 0 & 0 & 4 \\
1                & 4 & 0 & 0 & 0 & 0 & 0 & 0 & 4 \\
2                & 0 & 4 & 0 & 0 & 0 & 0 & 0 & 4 \\
3                & 0 & 3 & 1 & 0 & 0 & 0 & 0 & 4 \\
4                & 0 & 0 & 4 & 0 & 0 & 0 & 0 & 4 \\
5                & 0 & 0 & 2 & 2 & 0 & 0 & 0 & 4 \\
6                & 0 & 0 & 0 & 3 & 0 & 0 & 0 & 3 \\
7                & 0 & 0 & 0 & 2 & 2 & 0 & 0 & 4 \\
8                & 0 & 0 & 0 & 0 & 4 & 0 & 0 & 4 \\
9                & 0 & 0 & 0 & 0 & 1 & 3 & 0 & 4 \\
10               & 0 & 0 & 0 & 0 & 0 & 4 & 0 & 4 \\
11               & 0 & 0 & 0 & 0 & 0 & 0 & 4 & 4 \\
12               & 0 & 0 & 0 & 0 & 0 & 0 & 3 & 3 \\ \hline
$\f_c^a (k)$     & 8 & 7 & 7 & 7 & 7 & 7 & 7 & 50
\end{tabular}
\end{center}
\caption{The two-dimensional distribution and its marginals for $a=50, b=13, c=7$}\label{distribution50137}
\end{figure}

By analyzing the structure of the marginal distributions we can arrive at closed formulae for
$M_k(a;b)$. For example, one can immediately verify that
  \[ \f_{a;b} (j) = n = \fl{ \frac a b } , \ j = 0, \dots, b-1 \qquad \qquad \text{ if } a \equiv 0 \mod b \ , \]
and
  \begin{equation} \f_{a;b} (j) = \left\{ \begin{array}{cl} \fl{ \frac a b } + 1 & \text{ if } j=0 , \\
                                                  \fl{ \frac a b } & \text{ if } j>0 \end{array} \right. \qquad \qquad \text{ if } a \equiv 1 \mod b \end{equation}
Thus for $ a \equiv 0, 1 \mod b $ we immediately obtain
  \[ M_k (a;b) = \frac 1 a \fl{ \frac a b } \sum_{ j=1 }^{b-1} j^k \ . \]
In general, the one-dimensional frequencies can be bounded as
  \begin{equation}\label{ibounds} \fl{ \frac a b } \leq f_{a;b} (j) < \cl{ \frac a b } \ . \end{equation}
A book-keeping device that will help us keep track of the difference between the frequency
and $\fl{ \frac a b }$ is the following.
\begin{definition}
$ I_{a;b} (j) = \f_{a;b} (j) - \fl{ \frac a b } $.
\end{definition}
Notice that by (\ref{ibounds}) we have $ I_{a;b} (j) = 0, 1 $ for all $a,b, j=0, \dots, b-1$.
Accordingly we rewrite the $k\th$ moment of the Dilcher-Girstmair distribution as follows.
  \begin{equation}\label{mri} M_k (a;b) = \frac 1 a \left( \fl{ \frac a b } \sum_{ j=1 }^{b-1} j^k + \sum_{ j=1 }^{b-1} j^k I_{a;b} (j) \right) . \end{equation}
The second sum allows for a finer analysis of these moments. A trivial example follows from the fact that $I_{a;b} (j) \geq 0$:
  \[ M_k (a;b) \geq \frac 1 a \fl{ \frac a b } \sum_{ j=1 }^{b-1} j^k \ . \]
This bound gets achieved, for example, when $ a \equiv 0, 1 \mod b $.
In the following section, we study $I_{a;b} (j)$ and its second moments.


\section{Some Formulae and bounds for $M_2$}\label{dilchersect}

Of special interest is $M_2(a;b)$, due to its relationship to the Dedekind sum $\s(a,b)$.
According to (\ref{mri})
  \begin{equation}\label{mri2} a \, M_2(a;b) = \fl{ \frac a b } \frac{ (b-1)b(2b-1) }{ 6 } + \sum_{ j=1 }^{b-1} j^2 I_{a;b} (j) \ . \end{equation}
One may think about this identity in terms of the Dilcher-Girstmair distribution model:
Among the $a$ integers in $[0,a)$, we have at least $\fl{ a/b }$ of them in each interval
  \[ \left[ \frac{ka}{b} , \frac{(k+1)a}{b} \right) , \qquad k = 0, \dots, b-1 \ . \]
These integers are represented in the first term
on the right-hand side of (\ref{mri2}). Suppose $a \equiv l \mod b$ where $ 0 < l < b $ (the
case $b|a$ is special and very easy to handle: $I_{a;b}(j) = 0$ for all $j$);
then there are $l-1$ integers ``left" which haven't been accounted for (note that the first
interval $ [ 0 , a/b ) $ contains $ \fl{a/b} + 1 $ integers). These $l-1$ integers are
represented in the second term on the right-hand side of (\ref{mri2}).
In fact, one can say more about them. Because they are uniformly distributed among the $b$
intervals, we obtain
  \begin{equation}\label{d2} D_2(a;b) = \sum_{ j=1 }^{b-1} j^2 I_{a;b}(j) = \sum_{ m=1 }^{l-1} \fl{ \frac{mb}{l} }^2 = l \, M_2 (l;b) \ , \qquad l \geq 2 . \end{equation}
Note that, in particular, $D_2(a;b)$ depends on $a$ only via $l \equiv a \mod b$.
In special cases of $a \equiv l \mod b$, we can obtain closed formulas for $D_2(a;b)$,
given in the following theorem.

\begin{theorem}\label{exd2} Let $l \equiv a \mod b$ then $D_2(a;b)$ is given by the following formulae:
\begin{center}
\begin{tabular}{c|ll}
$l$ & $D_2(a;b)$ \\ \hline
0,1 & 0 \\
2   & $ \fl{ \frac b 2 }^2 $ \\
3   & $ 5 \fl{ \frac b 3 }^2 $ & if $ b \equiv 0, 1 \mod 3 $ \\
    & $ 5 \fl{ \frac b 3 }^2 + 4 \fl{ \frac b 3 } + 1 $ & if $ b \equiv 2 \mod 3 $ \\
4   & $ 14 \fl{ \frac b 4 }^2 $ & if $ b \equiv 0, 1 \mod 4 $ \\
    & $ 14 \fl{ \frac b 4 }^2 + 10 \fl{ \frac b 4 } + 2 $ & if $ b \equiv 2 \mod 4 $ \\
    & $ 14 \fl{ \frac b 4 }^2 + 16 \fl{ \frac b 4 } + 5 $ & if $ b \equiv 3 \mod 4 $ \\
5   & $ 30 \fl{ \frac b 5 }^2 $ & if $ b \equiv 0, 1 \mod 5 $ \\
    & $ 30 \fl{ \frac b 5 }^2 + 14 \fl{ \frac b 5 } + 2 $ & if $ b \equiv 2 \mod 5 $ \\
    & $ 30 \fl{ \frac b 5 }^2 + 26 \fl{ \frac b 5 } + 6 $ & if $ b \equiv 3 \mod 5 $ \\
    & $ 30 \fl{ \frac b 5 }^2 + 40 \fl{ \frac b 5 } + 14 $ & if $ b \equiv 4 \mod 5 $ \\
\end{tabular}
\end{center}
\end{theorem}

\emph{Proof.} This follows directly from (\ref{d2}).
\hfill {} $\Box$

In general, one can use (\ref{d2}) to obtain inequalities for $D_2(a;b)$ and hence for $M_2(a;b)$.
To this extend, we use the fact that
  \[ \left( m \fl{ \frac b l } \right)^2 \leq \fl{ \frac{mb}{l} }^2 \leq \fl { \left( \frac{mb}{l} \right)^2 } \ , \qquad l \geq 1 \ , \]
which implies the following bounds for $l \geq 2$.
  \[ \fl{ \frac b l }^2 \sum_{ m=1 }^{ l-1 } m^2 = \fl{ \frac b l }^2 \frac{ (l-1) l (2l-1) }{ 6 } \leq D_2 (a;b) \leq \fl { \left( \frac{b}{l} \right)^2 \frac{ (l-1) l (2l-1) }{ 6 } } = \fl { \sum_{ m=1 }^{ l-1 } \left( \frac{mb}{l} \right)^2 } \ . \]
Accordingly, we have the following:
\begin{theorem} For $a \equiv l \mod b$
  \begin{equation} M_2 (a;b) = \fl{ \frac a b } \frac{ (b-1)b(2b-1) }{ 6 } \qquad \text{ if } l=0,1, \end{equation}
and for $l \geq 2$
  \begin{equation}\label{flb1} M_2 (a;b) \geq \fl{ \frac a b } \frac{ (b-1)b(2b-1) }{ 6 } + \fl{ \frac b l }^2 \frac{ (l-1) l (2l-1) }{ 6 } \ . \end{equation}
and
  \begin{equation}\label{fub} M_2 (a;b) \leq \fl{ \frac a b } \frac{ (b-1)b(2b-1) }{ 6 } + \fl { \frac{ b^2 (l-1) (2l-1) }{ 6 l } } \ . \end{equation}
\end{theorem}

Naturally, (\ref{d2}) can be refined further to give even better bounds. We illustrate one further step here.
Suppose as before that $a \equiv l \mod b$ where $ 1 < l < b $, and moreover that $b = \fl{b/l} l + k$, so that
$0 \leq k \leq l-1$. According to (\ref{d2}) this gives
  \begin{align*} D_2(a;b) &= \sum_{ m=1 }^{l-1} \fl{ \frac{mb}{l} }^2
                           = \sum_{ m=1 }^{l-1} \fl{ \frac{m(\fl{b/l} l + k)}{l} }^2
                       = \sum_{ m=1 }^{l-1} \left( m \fl{ \frac b l } + \fl{ \frac{mk}{l} } \right)^2 \\
                          &= \fl{ \frac b l }^2 \frac{ (l-1) l (2l-1) }{ 6 } + 2 \fl{ \frac b l } \sum_{ m=1 }^{l-1} m \fl{ \frac{mk}{l} } + \sum_{ m=1 }^{l-1} \fl{ \frac{mk}{l} }^2 \ . \end{align*}

\begin{lemma} For all $l \geq 2$
  \[ D_2(a;b) = \fl{ \frac b l }^2 \frac{ (l-1) l (2l-1) }{ 6 } \qquad \text{ if } k = 0, 1 , \]
  \begin{equation}\label{flb2} D_2(a;b) = \fl{ \frac b l }^2 \frac{ (l-1) l (2l-1) }{ 6 } + \left\{ \begin{array}{ll} \fl{ \frac b l } \left( (l-1) l - \frac 1 4 l (l-2) \right) + \frac l 2 & \text{ if $l$ is even, } \\
                                                                                            \fl{ \frac b l } (l-1) l - \fl{ \frac l 2 } \left( 1 + \fl{ \frac l 2 } \right) + \fl{ \frac l 2 } & \text{ if $l$ is odd. } \end{array} \right. \ . \end{equation}
\end{lemma}
For example, if $l=4, k=2$
  \[ D_2(4,b) = \fl{ \frac b 4 }^2 \frac{ 3 \cdot 4 \cdot 7 }{ 6 } + 10 \fl{ \frac b 4 } + 2 \ , \]
and if $l=5, k=2$
  \[ D_2(5,b) = \fl{ \frac b 5 }^2 \frac{ 4 \cdot 5 \cdot 9 }{ 6 } + 14 \fl{ \frac b 5 } + 2 \ , \]
as stated in Theorem \ref{exd2}.
Finally, since $ \fl{ \frac{ mk }{ l } } \geq \fl{ \frac{ m2 }{ l } } $ for all $k \geq 2$, the above
formula of $D_2(a;b)$ for $k=2$ is a lower bound.

Similar bounds can be derived ``classically" by applying Dedekind's famous reciprocity law:
\begin{theorem}[Dedekind]\label{dedrec} If $a$ and $b$ are relatively prime then
  \begin{equation} \s (a,b) + \s (b,a) = - \frac{1}{4} + \frac{1}{12} \left( \frac{ a }{ b } + \frac{ 1 }{ ab } + \frac{ b }{ a } \right) \ . \end{equation}
\end{theorem}
Denote the rational function appearing in Theorem \ref{dedrec} by
  \begin{equation} R(a,b) = - \frac 1 4 + \frac 1 {12} \left( \frac a b + \frac 1 {ab} + \frac b a \right) \ . \end{equation}
Then we obtain for $a \equiv l \mod b$, where $a$ and $b$ are relatively prime and $ 1 < l < b $,
  \[ \s(b,a) = R(a,b) - \s(a,b) = R(a,b) - \s(l,b) = R(a,b) - R(b,l) + \s(b,l) \ . \]
It is well known (and a straightforward exercise) that
  \begin{equation} \left| \s(b,l) \right| \leq \s(1,l) = \frac l {12} - \frac 1 4 + \frac 1 {6l} \ , \end{equation}
which gives the following bounds:
  \begin{equation} R(a,b) - R(b,l) - \s(1,l) \leq \s(b,a) \leq R(a,b) - R(b,l) + \s(1,l) \ . \end{equation}
These inequalities, in turn, can be transformed into inequalities for $M_2$ via (\ref{M2}), to obtain:

\begin{theorem} Lower and upper bounds for $M_2$ are:
  \begin{equation}\label{rlb} M_2(a;b) \geq \frac{ (b^2 + 1) (a-1) (2a-1) }{6a^2} - \frac{ (a-1)b }{ 2a } - \frac{ 2b }{a} \left( R(a,b) - R(b,l) + \s(1,l) \right) . \end{equation}
  \begin{equation}\label{rub} M_2(a;b) \leq \frac{ (b^2 + 1) (a-1) (2a-1) }{6a^2} - \frac{ (a-1)b }{ 2a } - \frac{ 2b }{a} \left( R(a,b) - R(b,l) - \s(1,l) \right) . \end{equation}
\end{theorem}
In the following table we give the exact values of $a M_2(a;b)$ and their lower bounds.
We denote by

flb1 the lower bound according to (\ref{flb1}), \\
flb2 the lower bound according to (\ref{flb2}), \\
rlb the lower bound according to (\ref{rlb}), \\
fub the upper bound according to (\ref{fub}), \\
rub the upper bound according to (\ref{rub}).

Note that we can compute rlb and rub only when $a$ and $b$ are relatively prime.

\begin{center}
\begin{tabular}{c|c|c|c|c|c|c|c}
$a$ & $b$ & exact & flb1 & flb2 & rlb & fub & rub \\ \hline
5 & 2 & 2 & 2 & 2 & 2 & 2 & 2 \\
5 & 3 & 6 & 6 & 6 & 6 & 6 & 6 \\
5 & 4 & 14 & 14 & 14 & 14 & 14 & 14 \\
6 & 2 & 3 & 3 & 3 & & 3 & \\
6 & 3 & 10 & 10 & 10 & & 10 & \\
6 & 4 & 18 & 18 & 18 & & 18 & \\
6 & 5 & 30 & 30 & 30 & 30 & 30 & 30 \\
7 & 2 & 3 & 3 & 3 & 3 & 3 & 3 \\
7 & 3 & 10 & 10 & 10 & 10 & 10 & 10 \\
7 & 4 & 19 & 19 & 19 & 19 & 19 & 19.9 \\
7 & 5 & 34 & 34 & 34 & 34 & 34 & 34 \\
7 & 6 & 55 & 55 & 55 & 55 & 55 & 55 \\
35 & 7 & 455 & 455 & 455 & & 455 & \\
39 & 7 & 490 & 469 & 481 & 486.5 & 497 & 490 \\
40 & 7 & 501 & 485 & 501 & 498.2 & 513 & 503.8 \\
41 & 7 & 510 & 510 & 510 & 510 & 529 & 517.8 \\
10 & 3 & 15 & 15 & 15 & 15 & 15 & 15 \\
11 & 3 & 16 & 16 & 16 & 16 & 16 & 16 \\
21 & 6 & 185 & 185 & 185 & & 185 & \\
20 & 6 & 174 & 174 & 174 & & 174 & \\
11 & 7 & 126 & 105 & 117 & 122.5 & 133 & 126 \\
10 & 9 & 204 & 204 & 204 & 204 & 204 & 204 \\
11 & 9 & 220 & 220 & 220 & 220 & 220 & 220 \\
12 & 9 & 249 & 249 & 249 & & 249 & \\
13 & 9 & 260 & 260 & 260 & 260 & 274 & 264.5 \\
14 & 9 & 288 & 234 & 250 & 280.8 & 301 & 288 \\
15 & 9 & 315 & 259 & 286 & & 327 & \\
16 & 9 & 328 & 295 & 328 & 322.9 & 354 & 335.7 \\
17 & 9 & 344 & 344 & 344 & 344 & 381 & 359.75 \\
24 & 10 & 648 & 626 & 648 & & 657
\end{tabular}
\end{center}


\section{Bounds for generalized Dedekind sums: the case $d=2$}\label{d=2bounds}

\subsection{Applications of the Cauchy-Schwartz inequality}\label{cauchy-schwsubsec}

In the present section we discuss some relationships between the $S$ and the $M$-functions.
By definition, if $b_1 = \dots = b_d$
  \[ S_d (a; b \ \1_d) = M_d (a;b) \ , \ d = 1, 2, \dots \]
Here $\1_d$ denotes the $d$-dimensional vector all of whose
components are 1. The Cauchy-Schwartz inequality yields
immediately, for $d=2$, the inequality
  \[ S_2 (a; b_1, b_2) \leq \left( M_2 (a; b_1) M_2 (a; b_2) \right)^{1/2} \ , \]
with equality if and only if $b_1 = b_2$.
Let
  \begin{equation} R_2 (a; b_1, b_2) = \frac{ S_2 (a; b_1, b_2) }{ \sqrt{ M_2 (a; b_1) M_2 (a; b_2) } } \ . \end{equation}
In the following table we give a few values of $ R_2 (a; b_1, b_2) $.
\[ \begin{array}{r|r|c|c}
   a & b_2 & R_2 (a; 2, b_2) & R_2 (a; 3, b_2) \\
   \hline
   11 &  7 & 0.9163 & 0.9799 \\
   21 &  5 & 0.9237 & 0.9721 \\
   18 & 11 & 0.9297 & 0.9729 \\
   73 & 39 & 0.9189 & 0.9695 \\
   99 & 33 & 0.9192 & 0.9707 \\
   \end{array} \]

%

It is interesting to observe in this table that all these $R_2(a;b,c)$-values
are close to 1, and that among these values $R_2(a;2,b_2) < R_2(a;3,b_2)$.
The question is whether this inequality is always true. A partial answer is
given in Lemma \ref{partialanswer} of Section \ref{geometricrealiztion}.
Empirical evaluations lead us to the following conjecture:
\begin{conjecture}\label{firstconjecture} For all $a, b, c \geq 3$,
  \begin{equation} R_2 (a;b,c) \geq R_2 (5;2,3) = \frac{ \sqrt 3 }{ 2 } \ . \end{equation}
\end{conjecture}

Notice that according to the previous definitions, $R_2(a;b,c)$ is the
cosine of the angle between the two vectors
  \[\v_{b} = \left(\fl{ \frac{b}{a} },  \fl{ \frac{ 2 b }{a} }, \dots, \fl{\frac{(a-1)b}{a}}
\right) \qquad \text{ and } \qquad \v_{c} = \left(\fl{ \frac{c}{a} },  \fl{ \frac{ 2 c }{a} }, \dots, \fl{\frac{(a-1)c}{a}}
\right) \ . \]
In section \ref{geometricrealiztion} we present the geometrical correspondence,
which is utilized to obtain further results.

\begin{figure}[htb]
\begin{center}
\includegraphics[totalheight=4in]{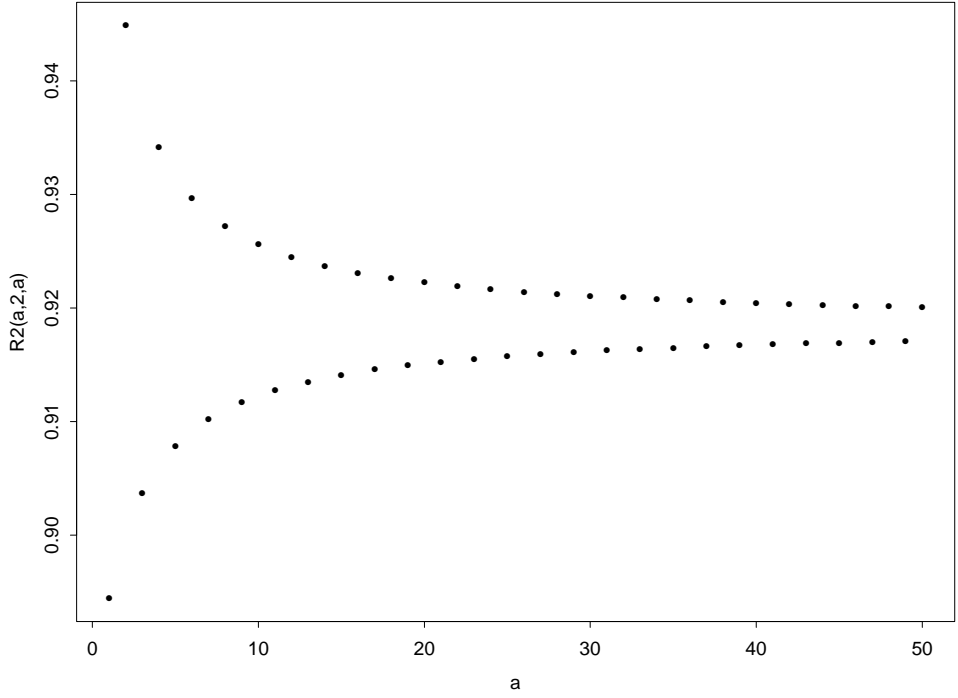}
\end{center}
\caption{$R_2(a;2,a)$ for $a=3, \dots, 50$}\label{plot}
\end{figure}

Exact formulae can be derived for $R_2 (a; 2, a), \ a \geq 3$. Indeed
  \begin{equation} S_2(a; 2, a) = \left\{ \begin{array}{cl} \frac{ a-1 }{ 2 } - \frac{ \fl{ a/2 } ( 1 + \fl{ a/2} ) }{ 2 a } & \text{ if } \fl{ \frac a 2 } < \frac a 2 , \text{ i.e., $a$ is odd, } \\
                                              \frac{ a-1 }{ 2 } - \frac 1 4 \left( \frac a 2 - 1 \right) & \text{ if } \fl{ \frac a 2 } = \frac a 2 , \text{ i.e., $a$ is even. } \end{array} \right. \end{equation}
Moreover
  \[ M_2 (a;2) = \frac 1 a \fl{ \frac a 2 } \ , \]
and
  \[ M_2 (a;a) = \frac{ 2 a^2 - 3 a + 1 }{ 6 } \ . \]
Accordingly
  \begin{equation} R_2(a; 2, a) = \left\{ \begin{array}{cl} R_2^* (a) \left( a - 1 - 1/a \fl{ a/2 } ( 1 + \fl{ a/2} ) \right) & \text{ if } \fl{ \frac a 2 } < \frac a 2 , \text{ i.e., $a$ is odd, } \\
                                              R_2^* (a) \left( a - 1 - 1/2 ( a/2 - 1 ) \right) & \text{ if } \fl{ \frac a 2 } = \frac a 2 , \text{ i.e., $a$ is even, } \end{array} \right. \end{equation}
where
  \begin{equation} R_2^* (a) = \frac{ \sqrt{ 6a } }{ 2 \sqrt{ \fl{ a/2 } ( 2 a^2 - 3 a + 1 ) } } \ . \end{equation}
A graph of $R_2(a;2,a)$ for $a=3, \dots, 50$ is given in Figure \ref{plot}. Notice that
$\lim_{a \to \infty} R_2(a;2,a) = \frac{ 3 \sqrt 6 }{ 8 }$.

We provide here a few auxiliary results. First, if $c=l+ia$ (i.e., $ a \equiv l \mod c $) then
  \begin{equation} S_2(a;b,c) = \sum_{m=0}^{a-1} \fl{ \frac{ bm }{ a } } \fl{ \frac{ (l+ia) m }{ a } } = i \sum_{m=1}^{a-1} m \fl{ \frac{ bm }{ a } } + S_2(a;b,l) \ . \end{equation}
Similarly,
  \[ a M_2(a;l+ia) = \sum_{m=0}^{a-1} \fl{ \frac{ (l+ia) m }{ a } }^2 = i^2 \frac{ (a-1) a (2a-1) }{ 6 } + 2i \sum_{m=1}^{a-1} m \fl{ \frac{ lm }{ a } } + a M_2(a,l) \ . \]
Accordingly,
  \begin{equation} R_2(a;b,l+ia) = \frac{ i \sum_{m=1}^{a-1} m \fl{ \frac{ bm }{ a } } + a S_2(a;b,l) }{ D_i } \end{equation}
where
  \[ D_i = i \left( a M_2(a;b) \frac{ (a-1) a (2a-1) }{ 6 } + \frac 2 i \sum_{m=1}^{a-1} m \fl{ \frac{ lm }{ a } } + \frac{ a }{ i^2 } M_2(a,l) \right)^{1/2} . \]
Thus,
  \[ \lim_{ c \to \infty } R_2(a;b,c) = \lim_{ i \to \infty} R_2(a;b,l+ia) = R_2(a;b,a) \ . \]
Now,
  \[ a S_2(a;ja,a) = \sum_{m=0}^{a-1} \fl{ \frac{ mja }{ a } } \fl{ \frac{ ma }{ a } } = j \sum_{m=0}^{a-1} m^2 = j \frac{ (a-1) a (2a-1) }{ 6 } \]
and
  \[ a M_2(a;ja) = j^2 \sum_{m=0}^{a-1} m^2 = j^2 \frac{ (a-1) a (2a-1) }{ 6 } \ , \]
whence
  \[ R_2(a;ja,a) = 1 \qquad \text{ for all } j \geq 1 \ . \]
We consider now $R_2(a;b,a)$ with $a \to \infty$. Let $a=jb$. For $j \geq 2$
  \[ R_2(jb;b,jb) = \frac{ \sum_{m=j}^{jb-1} m \fl{ \frac m j } }{ \left( \sum_{m=j}^{jb-1} \fl{ \frac m j }^2 \frac{ (a-1) a (2a-1) }{ 6 } \right)^{1/2} } \]
  \[ \sum_{m=j}^{jb-1} m \fl{ \frac m j } = \sum_{l=1}^{b-1} l \sum_{m=lj}^{(l+1)j-1} m = \frac{ bj (b-1) (4bj+j-3) }{ 12 } \]
  \[ \sum_{m=j}^{jb-1} \fl{ \frac m j }^2 = j \sum_{l=1}^{b-1} l^2 = j \frac{ (b-1) b (2b-1) }{ 6 } \ . \]
Thus
  \begin{equation} R_2^* (b) = \lim_{j \to \infty} R_2(jb;b,jb) = \lim_{j \to \infty} \frac{ (b-1) \left( (4b+1) j - 3 \right) }{ 2 \left( (b-1) (2b-1) (jb-1) (2jb-1) \right)^{1/2} } = \frac{ \sqrt 2 (b-1) (4b+1) }{ 4b \sqrt{ (b-1) (2b-1) } } \ . \end{equation}
Some values of the limit are given in the table below.
\begin{figure}[htb]
\begin{center}
\begin{tabular}{c|l}
$b$ & $R_2^* (b)$ \\ \hline
2 & 0.918558 \\
3 & 0.96896 \\
4 & 0.9836
\end{tabular}
\end{center}
\caption{Some values of $R_2^* (b) = \lim_{j \to \infty} R_2(jb;b,jb)$}
\end{figure}


\subsection{ A Geometric correspondence}\label{geometricrealiztion}\label{geomsubsec}

We have seen that the Dedekind-like sums $S_2(a;b,c)$ and $M_2(a;b)$ can be considered as
inner products in $\R^{a-1}$.  Thus, for a given integer $a \geq 3$, we construct a polyhedral cone
$\C_a \subset \R^{a-1}$ that is defined by the positive real span of the vectors
  \[ \v_b = \left( \fl{ \frac b a } , \fl{ \frac {2b} a } , \dots, \fl{ \frac {(a-1)b} a } \right),  \ 1 \leq b < \infty. \]
As stated in the introduction, the significance of these vectors is that
$R_2(a;b,c)$ is the cosine of the angle between the two vectors $\v_b$ and $\v_c$.
The observation that $R_2(a;b,c)$ is close to $1$ is captured geometrically by the statement
that this cone $\C_a$ is thin in the angular metric.

Notice that $\v_1 = \0 = (0, \dots, 0)$ and $\v_a = ( 1, 2, \dots, a-1 )$.
For $b=2, 3, \dots, a-1$, the vectors $\v_b$ are
  \[ \v_2 = (0, \dots, 0, 1, \dots, 1) , \v_3 = (0, \dots, 0, 1, \dots, 1, 2, \dots, 2) , \dots, \v_{a-1} = (0, 1, \dots, a-2), \]
where each vector $\v_j$ has almost equally distributed values for the integers the comprise its entries.

For $b > a$, we write $ b = ka + l, \ k > 0 , \ 0 \leq l < a $, and it follows from our notation
that  $\v_b = k \v_a + \v_l$.  Thus the cone $\C_a$ is in fact the positive real span of
only the $a-1$ vectors $\v_b$ with $b = 2, 3, \dots , a$.

Since $\v_1 = 0$, we have $\v_{ka} = \v_{ka+1}$. Moreover, if $P_l \ (2 \leq l \leq a-1)$
denotes the 2-dimensional plane containing the vectors $\v_l$ and $\v_a$, then all the vectors
$\v_{ka+l}, \ k = 0, 1, \dots$ belong to $P_l$.
Notice that for different values of $l$, say $l$ and $l' \not= l$, $P_l$ and $P_{l'}$
are two different planes which have the ray
  \[\left\{ r\, \v_a : \ r \geq 0 \right\} \]
in common (see Lemma \ref{firstlemma}).
Throughout this section, the denominators in the vector components of all the vectors
$\v_m$ are always the integer $a$.
\vspace{12pt}
\begin{figure}[ht]
\begin{center}
\includegraphics{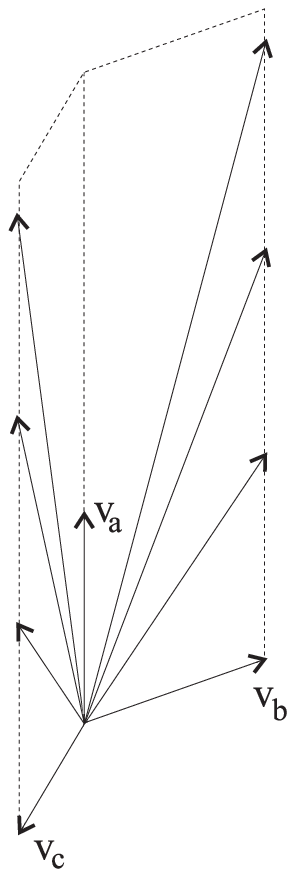}
\end{center}
\vspace{-12pt}
\caption{$\v_a, \ \v_b + k \v_a (0 \leq k \leq 3), \ \v_c + j \v_a (0 \leq j \leq 3)$}
\end{figure}


%

The vectors $\v_2, \dots, \v_a$ are not always linearly independent.
One can easily check that if $a=3,4,6$ then these vectors are linearly
independent, and when $a=5,7,8,9,\dots$ they are not.
However, one can prove the following:
\begin{lemma}\label{firstlemma}
For each $a \geq 4$, the vectors $\v_a, \v_l, \v_{l'}$,
where $1 < l, l' < a$ and $l \not= l'$, are linearly independent.
\end{lemma}
\emph{Proof.} For each $a \geq 4$, $\v_a = (1, 2, \dots, a-1)$, while
the first component of both $\v_l$ and $\v_{l'}$ is zero; thus $\v_a$
is linearly independent of $\{ \v_l, \v_{l'} \}$. Moreover, $\v_l$ and
$\v_{l'}$ do not lie on the same ray.
\hfill {} $\Box$

\begin{lemma}\label{secondlemma}
For each $a \geq 3$, if $b_1 = k_1 a + l$ and $b_2 = k_2 a + l$,
where $1 \leq k_1 < k_2$, then
  \[ R_2 (a;l,a) < R_2 (a;l,b_2) < R_2 (a;l,b_1) \]
for all $1 < l < a$.
\end{lemma}
\emph{Proof.} The vectors $\v_{b_1}$ and $\v_{b_2}$ lie in $P_l$.
Moreover, $\v_{b_1} = k_1 \v_a + \v_l$ and $\v_{b_2} = k_2 \v_a + \v_l$.
Hence
  \[ \an{ \v_l, \v_a } > \an{ \v_{b_1}, \v_a } > \an{ \v_{b_1}, \v_a } \ . \]
But $ \an{ \v_{b}, \v_l }  = \an{ \v_a, \v_l } - \an{ \v_a, \v_b } $,
and the inequalites in the statement follow by taking cosines.
\hfill {} $\Box$

Notice that due to the monotonicity stated in the last lemma,
  \[ \lim_{ k \to \infty } R_2 ( a; l, ka + l ) = R_2 ( a; l, a ) \ . \]
It is interesting to notice that in the case of $a=3$, all vectors $\v_b$
are between $\v_2 = (0,1)$ and $\v_3 = (1,2)$. The cosine of the angle
between these two vectors is $2 / \sqrt 5$.

\begin{lemma} \begin{enumerate}[{\rm (i)}]
\item\label{firstcase} If $a = 3, 5$ then $R_2 (a;2,b) < R_2 (a;3,b)$ for all $b \geq 3$.
\item\label{secondcase} If $a = 4$ then $R_2 (4;2,b) < R_2 (4;3,b)$ for all $b \not= 6$.
  If $b = 6$ then $R_2 (4;2,6) = 0.9708$ and $R_2 (4;3,6) = 0.9647$.
\end{enumerate}
\end{lemma}
\emph{Proof.} (\ref{firstcase})
The case $a=3$ follows immediately from Lemma \ref{secondlemma}.
For $a=5$ we have
\begin{align*} \v_2 = (0,0,1,1) , &\qquad | \v_2 | = \sqrt 2 \\
               \v_3 = (0,1,1,2) , &\qquad | \v_3 | = \sqrt 6 \\
               \v_4 = (0,1,2,3) , &\qquad | \v_4 | = \sqrt {14} \\
               \v_5 = (1,2,3,4) , &\qquad | \v_5 | = \sqrt {30} \ . \\
\end{align*}
Let $\sp \cdot \cdot$ denote the inner product of two vectors. For any
$b = 5k + l , \ k = 1, 2, \dots, \ l = 2, 3, 4$,
  \[ R_2 ( 5; 2, 5k + l ) < R_2 ( 5; 3, 5k + l ) \]
if and only if
  \[ \frac{ 1 }{ \sqrt 2 } \sp{ \v_2 }{ \v_{ 5k + l } } < \frac{ 1 }{ \sqrt 6 } \sp{ \v_3 }{ \v_{ 5k + l } } \ . \]
This is equivalent to
  \[ \sqrt 6 k \sp{ \v_2 }{ \v_5 } + \sqrt 6 \sp{ \v_2 }{ \v_l } < \sqrt 2 k \sp{ \v_3 }{ \v_5 } + \sqrt 2 \sp{ \v_3 }{ \v_l } \ , \]
or
  \[ k \left( 13 \sqrt 2 - 7 \sqrt 6 \right) > \sqrt 6 \sp{ \v_2 }{ \v_l } - \sqrt 2 \sp{ \v_3 }{ \v_l } \ . \]
Thus, for $l=2$ the inequality is true for all $k > 0.53$; for $l=3$ or 4 it is true for all
$k \geq 0$. Notice that for $l=2$ and $k=0$, we get $b=2$.

(\ref{secondcase}) For $a=4$, if $b=4k+2$ then the inequality is true for
all $k > 1.72$. For this reason, the inequality between $R_2(4;2,6)$ and
$R_2(4;3,6)$ is reversed. If $b=4k+3$ the inequality is true for all $k \geq 0$.
\hfill {} $\Box$

Empirical evidence suggests that $R_2 (a; 2, b) < R_2 (a; 3, b)$ for
all $a \geq 6$ and $b \geq 3$. We do not give a formal proof.

\begin{lemma} For all $a \geq 3$, $R_2 (a;2,a) < R_2 (a;3,a)$.
\end{lemma}
\emph{Proof.} If $a=3$ then $R_2 (3;3,3) > R_2 (3;2,3)$.
For all $ a \geq 4 $, we have to show that $B(a) > A(a)$, where
  \[ A(a) = R_2 (a;2,a) | \v_a | \qquad \text{ and } \qquad B(a) = R_2 (a;3,a) | \v_a | \ . \]
The formulas for $A(a) , \ a \equiv 0, 1 \mod 2$ and
$B(a) , \ a \equiv 0, 1, 2 \mod 3$ are given in Figure \ref{mods}.

\begin{figure}[htb]
\begin{center}

\begin{tabular}{c|l}
$a \mod 2$ & $A$ \\ \hline
0 & $\displaystyle \frac{\sqrt{2a}}{8} \, (3a-2) $ \\
1 & $\displaystyle \frac{ 1 }{ 2 \sqrt{ \fl{ a/2 } } } \left( a + \fl{ a/2 } \right) \left( a - 1 - \fl{ a/2 } \right) $ \\
\end{tabular}

\vspace{12pt}

\begin{tabular}{c|l}
$a \mod 3$ & $B$ \\ \hline
0 & $\displaystyle \frac{ \sqrt{ 15 a } }{ 90 } \, ( 13a - 9 ) $ \\
1 & $\displaystyle \frac{ a (a-1) - 3/2 \fl{ a/3 } - 5/2 \fl{ a/3 }^2 }{ 2 \sqrt{ a - 1 - 7/2 \fl{ a/3 } } } $ \\
2 & $\displaystyle \frac{ a (a-1) - 1 - 7/2 \fl{ a/3 } - 5/2 \fl{ a/3 }^2 }{ \sqrt{ 4a + 7 - 7 \fl{ a/3 } } } $ \\
\end{tabular}

\end{center}
\caption{$A(a) , \ a \equiv 0, 1 \mod 2$ and $B(a) , \ a \equiv 0, 1, 2 \mod 3$}\label{mods}
\end{figure}

One can easily check in all six cases that $B(a) > A(a)$.
\hfill {} $\Box$

\begin{lemma}\label{partialanswer}
For each $a \geq 3$ and each $b, c \geq 2$
  \[ R_2 (a;b,c) \geq \min \left( R_2 (a;l',a) , R_2 (a;l,a) , R_2 (a;l',l) \right) \]
where $b \equiv l \mod a , \ c \equiv l' \mod a$.
\end{lemma}
\emph{Proof.}
If $l, l' \equiv 0, 1 \mod a$, both $\v_b$ and $\v_c$ are on the ray $R_a$,
and $R_2 (a;b,c) = 1$.

If $l = l' = 2, \dots, a-1$ then $\v_b$ and $\v_c$ belong to $P_l$ and
$R_2 (a;b,c) \geq R_2 (a;l,a)$.

Finally, if $l \not= l'$ one establishes the inequality by comparing the
arcs on the unit sphere corresponding to the angles. These are the arcs
between the points on the sphere on the rays generated by $\v_a$, $\v_l$,
$\v_{l'}$, $\v_b$, and $\v_c$.
\hfill {} $\Box$

\begin{figure}[htb]
\begin{center}
\includegraphics[angle=270,totalheight=4in]{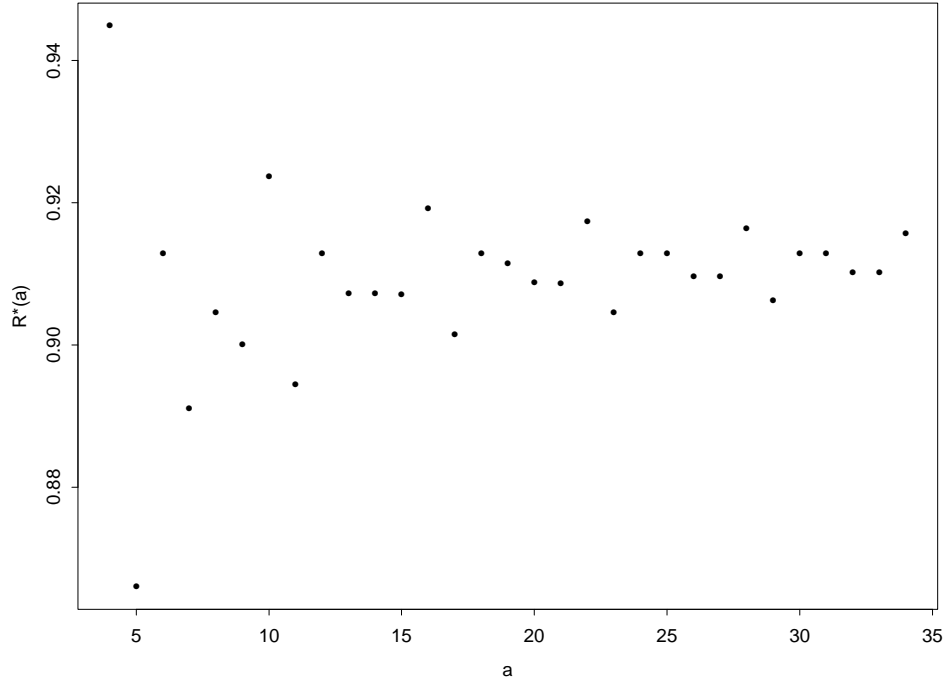}
\end{center}
\caption{$R_2^* (a)$ for $a=3, \dots, 35$}\label{r2starplot}
\end{figure}

To prove Conjecture \ref{firstconjecture} it suffices to show that
  \[ \min_{ a \geq 3, \ b,c \geq 2 } R_2 (a;b,c) \geq R_2 (5;2,3) \ . \]
Let $\displaystyle R_2^* (a) = \min_{ b,c \geq 2 } R_2 (a;b,c)$. According to Lemma
\ref{partialanswer}, for each $a \geq 3$
  \[ R_2^* (a) = \min_{ 1 < l, l' < a } \min \left( R_2 (a;l',a) , R_2 (a;l,a) , R_2 (a;l',l) \right) \ , \]
whence
  \[ \min_{ a \geq 3, \ b,c \geq 2 } R_2 (a;b,c) = \min_{ a \geq 3 } R_2^* (a) \ . \]
In Figure \ref{r2starplot} we present a plot of $R_2^* (a)$ for $a=3, \dots, 35$.
We see that in this range, $R_2 (5;2,3)$ is the minimum, lending further credence to conjecture 1.

\section{Upper bounds for generalized Dedekind sums: higher dimensions}\label{boundshigherdim}

\subsection{Probability models}

We introduce now a probability space and random variables, whose (mixed)
moments yield the $S$- and $M$-functions. Let $\D_d^{(a)}$ be a $d$-dimensional
discrete sample space, consisting of $a^d$ points, that is,
  \[ \D_d^{(a)} = \left\{ ( m_1 , \dots , m_d ) : \ m_j = 0, \dots , a-1 , \ j=1, \dots, d \right\} . \]
A point in $\D_d^{(a)}$ is a $d$-dimensional vector $\m = ( m_1 , \dots , m_d ) $.
Consider the probability function on $\D_d^{(a)}$:
  \begin{equation}\label{probfunction} P ( \m ) = \left\{ \begin{array}{cl} 1/a & \text{ if } \m = j \1_d , \ j=0, \dots, a-1 \\
                                                                                0 & \text{ otherwise, } \end{array} \right. \end{equation}
where $\1_d = (1, 1, 1, \dots, 1)$.
This probability function is concentrated on the main diagonal points of $\D_d^{(a)}$.
Define the random variables
  \begin{equation}\label{randomvariable} X_i^{(a)} ( \m ; \b ) = \left\{ \begin{array}{cl} \fl{ \frac{ m b_i }{ a } } & \text{ if } \m = m \1_d , \ m=0, \dots, a-1 , \\
                                                                                                                        0 & \text{ otherwise, } \end{array} \right. \end{equation}
where $\b = ( b_1, \dots, b_d )$.
It follows immediately that
  \begin{equation}\label{randomvariablecons} S_d (a; \b) = E_P \left\{ \prod_{i=1}^d X_i^{(a)} ( \m ; \b ) \right\} , \end{equation}
where $ E_P \{ \ \} $ denotes the expected value of the term in braces with respect to the
probability function $P$. Moreover,
  \begin{equation}\label{Mkexpected} M_k (a;b_i) = E_P \left\{ \left( X_i^{(a)} ( \m ; \b ) \right)^k \right\} . \end{equation}
Notice that $M_k (a;b_i)$ is the $k$-th order moment of $X_i^{(a)} ( \m ; \b )$.
The Dilcher-Girstmair presentation of the $S$- and $M$-functions can be described
as moments of the random variables
  \begin{equation}\label{Jdescription} J_i^{(a)} ( \m ; \b ) = \sum_{j=0}^{b_i - 1} j \, I \left\{ \m : \ \frac{ ja }{ b_i } \leq m_i < \frac{ (j+1) a }{ b_i } \right\} , \ i=1, \dots, d , \end{equation}
where $I \{ \m : \dots \}$ is the indicator function.
According to this definition,
  \begin{equation}\label{Sdexpected} S_d ( a; \b ) = E_P \left\{ \prod_{i=1}^d J_i^{(a)} ( \m ; \b ) \right\} \end{equation}
and
  \begin{equation}\label{MkexpectedJ} M_k (a;b_i) = E_P \left\{ \left( J_i^{(a)} ( \m ; \b ) \right)^k \right\} . \end{equation}


\subsection{Upper bounds for $S_d ( a; \b )$}

In the present section we use the random variables $X_i^{(a)} ( \m ; \b )$.
Since $a$ and $\b$ are fixed, we will simplify the notation to calling the
random variables $X_1, \dots, X_d$.
Repeated application of the Cauchy-Schwartz inequality yields bounds in terms
of the one-dimen\-sional moments $M$. For example, for $d=2$ we obtain
  \[ E \left\{ X_1 X_2 \right\} \leq \left( E \left\{ X_1^2 \right\} E \left\{ X_2^2 \right\}  \right)^{1/2} , \]
and thus
  \[ S_2 (a; b_1, b_2) \leq \left( M_2(a,b_1) M_2 (a,b_2) \right)^{1/2} . \]
For $d=3$ we get
  \[ E \left\{ X_1 X_2 X_3 \right\} \leq \left( E \left\{ X_1^2 \right\} E \left\{ X_2^2 X_3^2 \right\}  \right)^{1/2} \leq \left( E \left\{ X_1^2 \right\} \left( E \left\{ X_2^4 \right\} E \left\{ X_3^4 \right\} \right)^{1/2} \right)^{1/2} \]
or
  \begin{equation}\label{S3bound} S_3 (a; b_1, b_2, b_3) \leq M_2^{1/2} (a,b_1) M_4^{1/4} (a,b_2) M_4^{1/4} (a,b_3) \ . \end{equation}
By taking the geometric mean of the cyclical permutations, we get the symmetric upper bound
  \[ S_3 (a; b_1, b_2, b_3) \leq \left( \prod_{j=1}^3 M_2 (a,b_j) M_4 (a,b_j) \right)^{1/6} \ . \]
For $d=4$ we similarly obtain
  \begin{equation}\label{S4bound} S_4 (a; b_1, b_2, b_3, b_4) \leq \left( \prod_{j=1}^4 M_4 (a,b_j) \right)^{1/4} \ . \end{equation}
For $d=5$ we start with
  \begin{align*} S_5 ( a; \b ) &= E_P \left\{ X_1 \cdots X_5 \right\} \\
                                 &\leq \left( E_P \left\{ X_1^2 X_2^2 \right\} \right)^{1/2} \left( E_P \left\{ X_3^2 X_4^2 X_5^2 \right\} \right)^{1/2} \\
                                 &\leq \left( M_4 (a; b_1) M_4 (a; b_2) \right)^{1/4} \left( \prod_{j=3}^5 M_4 (a; b_j) M_8 (a; b_j) \right)^{1/12}
  \end{align*}
Symmetrizing this upper bound by taking the geometric mean of the $\binom 5 2 $ different
bounds obtained by different selections of pairs and triplets gives
  \begin{equation}\label{S5bound} S_5 ( a; \b ) \leq \left( \prod_{j=1}^5 M_4^3 (a; b_j) M_8 (a; b_j)  \right)^{1/20} . \end{equation}
From the upper bound for $S_3$ we immediately obtain
  \begin{equation}\label{S6bound} S_6 ( a; \b ) \leq \left( \prod_{j=1}^6 M_4 (a; b_j) M_8 (a; b_j)  \right)^{1/12} . \end{equation}
Generally, if $d=2k, \ k=1, 2, \dots$, we have
  \begin{equation}\label{S2kbound} S_{2k} ( a; \b ) \leq \left( E_P \left\{ X_1^2 \cdots X_k^2 \right\}  \right)^{1/2} \left( E_P \left\{ X_{k+1}^2 \cdots X_{2k}^2 \right\}  \right)^{1/2} , \end{equation}
from which we get, by utilizing previous results, symmetric upper bounds.
For example,
  \begin{align} S_{8} ( a; \b ) &\leq \left( M_8 (a; b_1) \cdots M_8 (a; b_4) \right)^{1/8} \left( M_8 (a; b_5) \cdots M_8 (a; b_8) \right)^{1/8} \label{S8bound} \\
                                    &= \left( M_8 (a; b_1) \cdots M_8 (a; b_8) \right)^{1/8} , \nonumber
  \end{align}
and
  \begin{equation}\label{S10bound} S_{10} ( a; \b ) \leq \left( \prod_{j=1}^{10} M_8^3 (a; b_j) M_{16} (a; b_j) \right)^{1/40} . \end{equation}
We can immediately prove by induction the following:

\begin{lemma}
  \begin{equation}\label{S2^kbound} S_{2^k} ( a; \b ) \leq \left( \prod_{j=1}^{ 2^k } M_{ 2^k } (a;b_j) \right)^{ 1/2^k } , \ k=1, 2, \dots \end{equation}
\end{lemma}
Similarly, for $k=0, 1, \dots$
  \begin{equation}\label{S32^kbound} S_{3 \cdot 2^k} ( a; \b ) \leq \left( \prod_{j=1}^{ 3 \cdot 2^k } M_{ 2^{k+1} } (a;b_j) M_{ 2^{k+2} } (a;b_j) \right)^{ 1/6 \cdot 2^k } \end{equation}
and
  \begin{equation}\label{S52^kbound} S_{5 \cdot 2^k} ( a; \b ) \leq \left( \prod_{j=1}^{ 5 \cdot 2^k } M_{ 2^{k+2} }^3 (a;b_j) M_{ 2^{k+3} } (a;b_j) \right)^{ 1/20 \cdot 2^k } . \end{equation}
If $d=2k+1$ one needs a two-stage process of first partitioning to
  \[ \left( E_P \left\{ \prod_{j=1}^k X_j^2 \right\}  \right)^{1/2} \left( E_P \left\{ \prod_{j=1}^{k+1} X_{k+j}^2 \right\}  \right)^{1/2} \]
and then symmetrizing.

Before concluding this section, we remark that the above upper bounds for the
$S$-functions are generally not unique. By different partitions one can obtain
different bounds. For example, in the case of $S_5$, one could start with
  \[ E_P \left\{ X_1 \cdots X_5 \right\} \leq \left( E_P \left\{ X_1^2 \right\}  \right)^{1/2} \left( E_P \left\{ X_{2}^2 \cdots X_{5}^2 \right\}  \right)^{1/2} = \left( M_2 (a;b_1) \right)^{1/2} \left( \prod_{j=2}^5 M_8 (a;b_j) \right)^{1/8}  . \]
After symmetrization we get
  \begin{equation}\label{S5boundalt} S_5 ( a; \b ) \leq \left( \prod_{j=1}^5 M_2 (a; b_j) M_8 (a; b_j) \right)^{1/10} . \end{equation}
The question is which upper bound should be used, (\ref{S5bound}) or (\ref{S5boundalt})?
For example, if $a=31$ and $\b = ( 3, 5, 7, 11, 13 ) $ then $S_5 ( a; \b ) = 1213.806$.
The upper bound given by (\ref{S5bound}) is $1321.321$, whereas that given by (\ref{S5boundalt})
is $1456.985$.
In the following table we present some exact values of $S_5 ( a; \b )$ and the two
bounds (\ref{S5bound}) and (\ref{S5boundalt}). We also show $R_5 ( a; \b )$, the
ratio of $S_5 ( a; \b )$ to the upper bound (\ref{S5bound}).

\begin{figure}[htb]
\begin{center}
\begin{tabular}{c|c|c|c|c|c}
$a$ & $\b$ & $S_5$ & bound (\ref{S5bound}) & bound (\ref{S5boundalt}) & $R_5$ \\ \hline
31 & (3, 5, 7, 11, 13) & 1213.806 & 1321.321 & 1456.985 & 0.9186 \\
21 & (5, 7, 9, 11, 13) & 4411.333 & 4668.719 & 5190.201 & 0.9449 \\
23 & (5, 9, 11, 13, 17) & 11429.74 & 12050.58 & 13385.72 & 0.9485 \\
27 & (5, 11, 13, 17, 21) & 28101.93 & 29617.94 & 33011.8 & 0.9488 \\
33 & (7, 11, 13, 19, 23) & 51943.76 & 54384.26 & 60525.59 & 0.9551
\end{tabular}
\end{center}
\caption{Some values and bounds of $S_5$}\label{S5table}
\end{figure}

It seems from Figure \ref{S5table} that the upper bound given by (\ref{S5bound})
is closer to the exact value of $S_5 ( a; \b )$ than (\ref{S5boundalt}). It is
the preferred upper bound. It is also interesting that, like in the case of
$R_2 (a; \b)$, all values of $R_5 ( a; \b )$ in Figure \ref{S5table} are
greater than $0.9186$.


\subsection{Relationships to upper bounds revisited}

We study now upper bounds to $S_d$ of the type given by (\ref{S3bound})--(\ref{S52^kbound}).
In particular, define
\begin{align} R_3 ( a; \b ) &= \frac{ S_3 (a; \b) }{ \left( \prod_{j=1}^3 M_2 (a,b_j) M_4 (a,b_j) \right)^{1/6} } \ , \\
              R_4 ( a; \b ) &= \frac{ S_4 (a; \b) }{ \left( \prod_{j=1}^4 M_4 (a,b_j) \right)^{1/4} } \ , \text{ and } \\
              R_5 ( a; \b ) &= \frac{ S_5 (a; \b) }{ \left( \prod_{j=1}^5 M_4^3 (a; b_j) M_8 (a; b_j)  \right)^{1/20} } \ .
\end{align}
A few values of $R_5 ( a; \b )$ are given in in Figure \ref{S5table}.
It seems that the minimal $R_5 ( a; \b )$ value is $R_5 ( 7; 2, 3, 4, 5, 6 ) = 0.8567$.
It is also interesting to observe that $R_5$, as shown in Figure \ref{S5table},
is generally above $0.9$, as in the case of $R_2$, despite the increase in
dimension from 2 to 5. We try to explain this phenomenon in probability terms.

As shown in (\ref{randomvariablecons}) and (\ref{Sdexpected}),
  \[ S_d ( a; \b ) = E_P \left\{ X_1 \cdots X_d \right\} = E_P \left\{ J_1 \cdots J_d \right\} \ . \]
Consider the case $d=2$. According to the law of iterated expectation \cite{bickeldoksum},
  \begin{equation}\label{iteratedexp} S_2 ( a; b_1, b_2 ) = E_P \left\{ J^{(a)} ( \m; b_1 ) E_P \left\{ J^{(a)} ( \m; b_2 ) \mid J^{(a)} ( \m; b_1 ) \right\} \right\} \ , \end{equation}
where the second term on the right-hand side is the \emph{conditional expectation} of
$J^{(a)} ( \m; b_2 )$, given $J^{(a)} ( \m; b_1 )$.
Notice that in the notation of Section \ref{2dimsect},
  \begin{equation}\label{S2rev} S_2 ( a; b, c ) = \frac 1 a \sum_{ j=0 }^{ b-1 } \sum_{ l=0 }^{ c-1 } j l f_{a;b,c} (j,l) = \frac 1 a \sum_{ j=0 }^{ b-1 } j f_{a;b} (j) \sum_{ l=0 }^{ c-1 } l \frac{ f_{a;b,c} (j,l) }{ f_{a;b} (j) } \ . \end{equation}
The key to understanding the phenomenon is that the joint frequencies $f_{a;b,c} (j,l)$
are distributed along the main diagonal, as illustrated in Figure \ref{distribution50137}.
In the special case that $b=c$,
  \[ f_{a;b,c} (j,l) = \left\{ \begin{array}{cl} f_{a;b} (j) = f_{a;c} (j) & \text{ if } j \not= l , \\
                                                                         0 & \text{ otherwise. } \end{array} \right. \]
In this case,
  \[ \sum_{ l=0 }^{ b-1 } l \frac{ f_{a;b,c} (j,l) }{ f_{a;b} (j) } = j \]
and
  \[ S_2 ( a; b, b ) = \frac 1 a \sum_{ j=0 }^{ b-1 } j^2 f_{a;b} (j) = M_2 ( a; b ) \ , \]
as expected.
When $b \not= c$ then $R_2$ is always smaller than 1, but might be quite close to it,
even when $b$ and $c$ are different. For example, $R_2 ( 50; 7, 13 ) = 0.9955$.


\subsection{One-dimensional moments relationships}

We present here some inequalities between $M_r (a; b)$ for fixed values of $a$ and $b$ and
for variable $r$.

\begin{theorem} $M_r (a; b)$ is log-convex in $r$.  That is,
  \begin{equation}  M_{ 2r } (a; b) M_{ 2r+2 } (a; b) - M_{ 2r+1 }^2 (a; b) > 0. \end{equation}
\end{theorem}

\emph{Proof.}
First, by Liapounov's inequality of moments \cite[p.~627]{kotzjohnson} we have
  \[ M_1 < M_2^{1/2} < M_3^{1/3} < \dots \]
By factoring $ \left\lfloor \dots \right\rfloor^{ 2r+1 } =
\left\lfloor \dots \right\rfloor^{ r } \left\lfloor \dots \right\rfloor^{ r+1 } $
we obtain the inequality
 \begin{equation} M_{ 2r+1 }^2 (a; b) < M_{ 2r } (a; b) M_{ 2r+2 } (a; b) \end{equation}
for all $r \geq 1$. That is, $M_r (a;b)$ is log-convex in $r$.
\hfill {} $\Box$


\section*{Appendix}

\begin{center}
{\bf Algorithm for Determining $f_{b}^a (i) , \ i=0, \dots, b-1$}
\end{center}
\small
\begin{verbatim}
STEP 0:
Set:
    (i) r = a/b;
    (ii) l = a-b*[r];
    (iii) f = 0(1,b);         # b-dimensional vector of zeros

STEP 1:

Compute:
    f [0] <- 1+[r];
    f [ i ] <- [(i+1)*r]- [i*r], i=1,..,(b-2);
    f [ b-1] <- a-1- [(b-1)*r];

STEP 2:

    IF((l=0) or (b=2)) { GOTO STEP 3 };
    ELSE {
        FOR ( i=1,...,b-2) {
            IF( [(i+1)*r] = r*(i+1) ) {
                f[i] <- f[i]-1;
                f[i+1] <- f[i+1]+1;
            }
        }

STEP 3:

PRINT f

END.
\end{verbatim}

\begin{center}
{\bf Algorithm for Determining $f_{bc}^a (i,j)$}
\end{center}
\begin{verbatim}
STEP 0:

Set:
    r1 <- a/b;
    r2 <- a/c;
    l1 <- a -b*[r1];
    l2 <- a - c*[r2];
    CT <- O((b+1),(c+1));    # matrix of zeros, of dimensions (b+1)*(c+1)
    t1 <- O(b,1);
    t2 <- O(1,c);

STEP 1:

Compute:
    CT[i,(c+1)] <- f_b[i-1], i=1,...,b;
    CT[(b+1), j] <- f_c[j-1], j=1,...,c;
    CT[(b+1),(c+1)] <- a;
    CT[1,1] <- min( CT[1,c+1], CT[b+1,1]);
    t2[1] <- t2[1] + CT[1,1];

STEP 2:

Compute:
    FOR (i=2,...,b) {
    CT[i,1] <- max(0, min(f_c[1]-t2[1], f_b[i-1])));
    t2[1] <- t2[1] + CT[i,1];
    }
    t1 <- CT[1:b,1];

STEP 3:

Compute:
    FOR (j=2,...,c) {
    CT[1,j] <- max(0,min(f_b[0]-t1[1], f_c[j-1]));
    t1[1] <- t1[1] + CT[1,j];
    t2[j] <- t2[j] + CT[1,j];
    }

STEP 4:

Compute:
    FOR (i=2,...,b) {
        FOR(j=2,...,c) {
        cty<- max(0, min(f_b[i-1]- t1[i], f_c[j-1]));
        ctx <- max(0, min(f_c[j-1]-t2[j], f_b[i-1]));
        CT[i,j] <- min(ctx,cty)
        t1[ i ] <- t1 [ i ]+  CT [ i, j];
        t2[ j ] <- t2 [ j ] + CT [i , j];
        }
    }

STEP 5:
Print CT

END.
\end{verbatim}
\normalsize


\bibliographystyle{plain}
\bibliography{bib}

\setlength{\parskip}{0cm}

\end{document}